\documentclass[10pt]{amsart}

\usepackage{amssymb}

\usepackage{enumerate}

\newcommand{\dis}{\rho(\e,\e')}

\newcommand{\N}{\ensuremath{\mathbb{N}}}

\newcommand{\e}{\varepsilon}

\newtheorem{thm}{Theorem}[section]

\newtheorem{cor}[thm]{Corollary}

\newenvironment{rk}{\refstepcounter{thm}\noindent%

{\bf Remark \arabic{section}.\arabic{thm}} \ }{

\smallskip

}

\newenvironment{pf}[1][]{\noindent {\it Proof #1} : }{\hbox{~}\qed \smallskip}

\title{Metrical characterization of super-reflexivity and linear type of Banach spaces}

\date{}

\author{Florent Baudier$^{\dag}$}

\begin{document}

\maketitle

\pagestyle{myheadings} \markright{Metrical characterization of
super-reflexivity and linear type of Banach spaces}

\begin{abstract}
We prove that a Banach space $X$ is not super-reflexive if and only if the
hyperbolic infinite tree embeds metrically into $X$. We improve one implication
of J.Bourgain's result who gave a metrical characterization of
super-reflexivity in Banach spaces in terms of uniforms embeddings of the
finite trees. A characterization of the linear type for Banach spaces is given
using the embedding of the infinite tree equipped with the metrics $d_{p}$
induced by the $\ell_{p}$ norms.

\end{abstract}

\makeatletter

\renewcommand{\@makefntext}[1]{#1}

\makeatother

\footnotetext{ Laboratoire de Math\'{e}matiques, UMR 6623

 Universit\'{e} de Franche-Comt\'{e},

 25030 Besan\c con, cedex - France \\ \indent $^{\dag}$florent.baudier@math.univ-fcomte.fr\\
 \indent \subjclass{2000 {\it Mathematics Subject Classification}.}{ (46B20)  (51F99)}}

\section{Introduction and Notation}

We fix some notation and recall basic results.\\

Let $(M,d)$ and $(N,\delta)$ be two metric spaces and an injective map $f:M\to
N$. Following \cite{MN}, we define the {\it distortion} of $f$ to be

$$ {\rm dist}(f):= \|f\|_{Lip}\|f^{-1}\|_{Lip}=\sup_{x\neq y \in
M}\frac{\delta(f(x),f(y))}{d(x,y)}.\sup_{x\neq y \in
M}\frac{d(x,y)}{\delta(f(x),f(y))}.$$ If $\mbox{dist}(f)$ is finite, we say
that $f$ is a metric embedding, or simply an embedding of $M$ into $N$.\\
And if there exists an embedding $f$ from $M$ into $N$, with dist$(f)\leq C$,
we use the notation $M \buildrel {C}\over
{\hookrightarrow} N$.\\

Denote $\Omega_{0}=\{\varnothing\}$, the root of the tree. Let
$\Omega_{n}=\{-1,1\}^{n}$, $T_{n}=\bigcup_{i=0}^{n}\Omega_{i}$ and
$T=\bigcup_{n=0}^{\infty}T_{n}$.  Thus
$T_{n}$ is the finite tree with $n$ levels and $T$ the infinite tree.\\
\indent For $\e$, $\e'\in T$, we note $\e\le\e'$ if $\e'$ is an extension of $\e$.\\
\indent Denote $\vert\e\vert$ the length of $\e$; i.e the numbers
of nodes of $\e$. We define the hyperbolic distance between $\e$
and $\e'$ by
$\rho(\e,\e')=\vert\e\vert+\vert\e'\vert-2\vert\delta\vert$, where
$\delta$ is the greatest common ancestor of $\e$ and $\e'$.
The metric on $T_{n}$, is the restriction of $\rho$.\\

For a Banach space $X$, we denote $B_{X}$ its closed unit ball, and $X^{*}$ its
dual
space.\\

\indent $T$ embeds isometrically into $\ell_{1}(\N)$ in a trivial way.
Actually, let $(e_{\e})_{\e\in T}$ be the canonical basis of $\ell_{1}(T)$ ($T$
is countable),
then the embedding is given by $\e\mapsto\sum_{s\le\e} e_{s}$.\\

\indent Aharoni proved in \cite{A} that
every separable metric space embeds into $c_{0}$, so $T$ does.\\

\indent The main result of this article is an improvement of Bourgain's
metrical characterization of super-reflexivity. Bourgain proved in \cite{B}
that $X$ is not super-reflexive if and only if the finite trees $T_{n}$
uniformly embed into $X$ (i.e with embedding constants independent of $n$).
Obviously if $T$ embeds into $X$ then the $T_{n}'s$ embed uniformly into $X$
and $X$ is not super-reflexive, but if $X$ is not super-reflexive we did not
know whether the infinite tree $T$ embeds
into $X$. In this paper, we prove that it is indeed the case : \\

\begin{thm}\label{theo1}
Let $X$ be a non super-reflexive Banach space, then $(T,\rho)$ embeds into $X$.
\end{thm}

\indent The proof of the direct part of Bourgain's Theorem essentially uses
James' characterization of super-reflexivity (see \cite{J}) and an
enumeration of the finite trees $T_{n}$. We recall James' Theorem :\\

\begin{thm}[James]
Let $0<\theta<1$ and $X$ a non super-reflexive Banach space, then :\\
$\forall\ n\in\N$, $\exists\ x_{1},x_{2},\dots,x_{n}\in B_{X}$, $\exists\
x^{*}_{1},x^{*}_{2},\dots ,x^{*}_{n}\in B_{X^{*}}$ s.t :

$$\begin{array}{cc}

 x^{*}_{k}(x_{j})=\theta & \forall k< j\\
 x^{*}_{k}(x_{j})=0 & \forall k\ge j\\
\end{array}$$

\end{thm}

\section{Metrical characterization of super-reflexivity}

The main obstruction to the  embedding of $T$ into any non-super-reflexive
Banach space $X$ is the finiteness of the sequences in James' characterization.
How, with a sequence of Bourgain's type embedding, can we construct a single
embedding
from $T$ into $X$ ?\\
\indent In \cite{R}, Ribe shows in particular, that
$\bigoplus_{2}l_{p_{n}}$ and $(\bigoplus_{2}l_{p_{n}})\bigoplus
l_{1}$ are uniformly homeomorphic, where $(p_{n})_{n}$ is a
sequence of numbers such that $p_{n}>1$, and $p_{n}$ tends to $1$.
But $T$ embeds into $l_{1}$, hence via the uniform homeomorphism
$T$ embeds into $\bigoplus_{2}l_{p_{n}}$. However $T$ does not
embed into any
$l_{p_{n}}$(they are super-reflexive).\\
\indent The problem solved in the next theorem, inspired  in part by Ribe's
proof, is to construct a subspace with a Schauder decomposition $\bigoplus
F_{n}$ where $T_{2^{n+1}}$ embeds into $F_{n}$ and to repast properly the
embeddings in order
to obtain the desired embedding.\\
\\

\begin{pf}[of Theorem 1.1]
Let $(\e_{i})_{i\ge 0}$, a sequence of  positive real numbers such that\\
$\prod_{i\ge 0}(1+\e_{i})\le 2$, and fix $0<\theta<1$. Let $k_{n}=2^{2^{n+1}+1}-1$.\\
First we construct inductively a sequence $(F_{n})_{n\ge 0}$ of subspaces of
$X$, which is a Schauder finite dimensional decomposition of a subspace of $X$
s.t the projection from $\bigoplus_{i=0}^{q} F_{i}$ onto $\bigoplus_{i=0}^{p}
F_{i}$, with kernel $\bigoplus_{i=p+1}^{q} F_{i}$ (with $p< q$) is of norm at
most $\prod_{i=p}^{q-1}(1+\e_{i})$, and sequences
$$x_{n,1},x_{n,2},\dots,x_{n,k_{n}}\in B_{F_{n}}$$
$$x^{*}_{n,1},x^{*}_{n,2},\dots,x^{*}_{n,k_{n}}\in B_{X^{*}}$$ s.t :\\

$$\begin{array}{cc}

 x^{*}_{n,k}(x_{n,j})=\theta & \forall k< j\\
 x^{*}_{n,k}(x_{n,j})=0 & \forall k\ge j.\\

\end{array}$$

Denote $\Phi_{n}$ : $T_{n}\to \{1,2,\dots,2^{n+1}-1\}$ the enumeration of
$T_{n}$ following the lexicographic order. It is an enumeration of $T_{n}$ such
that any pair of segments in $T_{n}$ starting at incomparable nodes (with
respect to the tree ordering $\le$) are mapped inside disjoint intervals. \\

\noindent Let $\Psi_{n}=\Phi_{2^{n+1}}$ and $\Gamma_{n}=T_{2^{n+1}}$.\\

$X$ is non super-reflexive, hence from James' Theorem :\\
$\exists\ x_{0,1},x_{0,2},\dots,x_{0,7}\in B_{X}$, $\exists\
x^{*}_{0,1},x^{*}_{0,2},\dots,x^{*}_{0,7}\in B_{X^{*}}$ s.t :

$$\begin{array}{cc}

 x^{*}_{0,k}(x_{0,j})=\theta & \forall k< j\\
 x^{*}_{0,k}(x_{0,j})=0 & \forall k\ge j.\\
\end{array}$$

\noindent $\Gamma_{0}=T_2$ embeds into $X$ via the embedding $
f_{0}(\e)=\sum_{s\le\e}
x_{0,\Psi_{0}(s)}$ (see \cite{B}).\\
Let $F_{0}={\rm Span}\{x_{0,1},\dots,x_{0,7}\}$, then ${\rm dim}(F_{0})<\infty$. \\

\noindent Suppose that $F_{0},\dots,F_{p}$, and
$$x_{p,1},x_{p,2},\dots,x_{p,k_{p}}\in B_{F_{p}}$$
$$x^{*}_{p,1},x^{*}_{p,2},\dots,x^{*}_{p,k_{p}}\in B_{X^*}$$ verifying the required
conditions, are constructed for all $p\le n$. \\

\noindent We apply Mazur's Lemma (see \cite{LT} page $4$) to the finite
dimensional subspace $\bigoplus_{i=0}^{n} F_{i}$ of $X$. Thus there exists
$Y_{n}\subset X$ such that ${\rm dim}(X/Y_{n})<\infty$ and :$$\Vert
x\Vert\le(1+\e_{n})\Vert x+y\Vert, \forall (x,y)\in\bigoplus_{i=0}^{n}
F_{i}\times Y_{n}$$.\\

\noindent But $Y_{n}$ is of finite codimension in $X$, hence is not
super-reflexive.\\
\noindent From James' Theorem and Hahn-Banach Theorem:\\
$$\exists\ x_{n+1,1},x_{n+1,2},\dots,x_{n+1,k_{n+1}}\in B_{Y_{n}},$$
$$\exists\ x^{*}_{n+1,1},x^{*}_{n+1,2},\dots,x^{*}_{n+1,k_{n+1}}\in B_{X^{*}},$$ s.t :

$$\begin{array}{cc}

 x^{*}_{n+1,k}(x_{n+1,j})=\theta & \forall k< j\\
 x^{*}_{n+1,k}(x_{n+1,j})=0 & \forall k\ge j.\\
\end{array}$$

\noindent $\Gamma_{n+1}$ embeds into $Y_{n}$ via the embedding $
f_{n+1}(\e)=\sum_{s\le\e} x_{n+1,\Psi_{n+1}(s)}$ .\\
\noindent Let $F_{n+1}={\rm Span}\{x_{n+1,j}\ ;\ 1\le j\le k_{n+1}\}$, then
${\rm dim}(F_{n+1})<\infty$, which achieves the induction.\\

Now define the following projections :\\

\noindent Let, $P_{n}$ the projection from
$\overline{\rm{Span}}(\bigcup_{i=0}^{\infty}F_{i})$ onto
$F_{0}\bigoplus\dots\bigoplus F_{n}$ with kernel
$\overline{\rm{Span}}(\bigcup_{i=n+1}^{\infty}F_{i})$.\\

\noindent It is easy to show that $\Vert
P_{n}\Vert\le\prod_{i=n}^{\infty}(1+\e_{i})\le
2$.\\

\noindent We denote now $\Pi_0=P_0$ and $\Pi_n=P_n-P_{n-1}$ for $n\geq 1$. We
have that $\|\Pi_n\|\leq 4$.\\

From Bourgain's construction, for all $n$ :
\begin{equation}
\frac{\theta}{3}\dis\le\Vert f_{n}(\e)-f_{n}(\e')\Vert\le\dis,
\end{equation}
where $f_{n}$ denotes the Bourgain's type embedding from $\Gamma_{n}$ in
$F_{n}$, i.e $ f_{n}(\e)=\sum_{s\le\e} x_{n,\Psi_{n}(s)}$.\\

Note that :\\
$$\forall\ n, \forall\ \e\in\Gamma_n\ \Vert f_n(\e)\Vert\le\vert\e\vert.$$

Now we define our embedding.\\

Let
$$\begin{array}{cccl}

    f : & T & \rightarrow & Y=\overline{\rm{Span}}(\bigcup_{i=0}^{\infty}F_{i})\subset X \\
        &   &   &  \\
        & \e & \mapsto & \lambda f_{n}(\e)+(1-\lambda)f_{n+1}(\e)\ ,\ {\rm if}\ 2^{n}\le\vert\e\vert\le 2^{n+1} \\
    \end{array}$$

where, $$\lambda=\frac{2^{n+1}-\vert\e\vert}{2^{n}}$$\\
And $$f(\varnothing)=0.$$\\
We will prove that :

\begin{equation}
\forall \e,\e'\in T\ \frac{\theta}{24}\dis\le\Vert f(\e)-f(\e')\Vert\le 9\dis.
\end{equation}

\begin{rk}\label{rk1}
We have $\frac{\theta}{24}\vert\e\vert\le\Vert f(\e)\Vert\le\vert\e\vert$.\\
\end{rk}

First of all, we show that $f$ is $9-$Lipschitz.

\noindent We can suppose that $0<\vert\e\vert\le\vert\e'\vert$ w.r.t remark \ref{rk1}.\\

\noindent If $\vert\e\vert\le\frac{1}{2}\vert\e'\vert$ then :

$$\dis\ge\vert\e'\vert -\vert\e \vert\ge\frac{\vert\e\vert +\vert\e' \vert}{3}$$

\noindent Hence, $$\Vert f(\e)-f(\e')\Vert\le 3\dis.$$

\noindent If $\frac{1}{2}\vert\e'\vert<\vert\e\vert\leq \vert\e'\vert$, we have
two different cases to consider.\\

\begin{enumerate}[1)]

\item if $2^{n}\le\vert\e\vert\le\vert\e'\vert<2^{n+1}$.\\
Then, let
$$\lambda=\frac{2^{n+1}-\vert\e\vert}{2^{n}}\ \ {\rm and}\ \ \lambda'=\frac{2^{n+1}-\vert\e'\vert}{2^{n}}.$$
$$\begin{array}{ccl}
        \Vert f(\e)-f(\e')\Vert & = & \Vert \lambda
        f_{n}(\e)-\lambda'f_{n}(\e')+(1-\lambda)f_{n+1}(\e)-(1-\lambda')f_{n+1}(\e')\Vert\\
                 & \le & \lambda\Vert f_{n}(\e)-f_{n}(\e')\Vert+\vert\lambda-\lambda'\vert(\Vert
         f_{n}(\e')\Vert+\Vert
         f_{n+1}(\e')\Vert)+(1-\lambda)\Vert
         f_{n+1}(\e)-f_{n+1}(\e')\Vert\\
          & \le &
          \dis+2\dis+2\dis\\
                & \le & 5\dis,
         \end{array}$$
because $\Vert f_n(\e')\Vert<2^{n+1}$, $\Vert f_{n+1}(\e')\Vert<2^{n+1}$ and,
$$\vert\lambda-\lambda'\vert=\frac{\vert\e'\vert-\vert\e\vert}{2^{n}}\le\frac{\dis}{2^{n}}.$$

\item if $2^{n}\le\vert\e\vert\le 2^{n+1}\le\vert\e'\vert<
2^{n+2}$.\\

Then, let
$$\lambda=\frac{2^{n+1}-\vert\e\vert}{2^{n}}\ \ {\rm and}\ \ \lambda'=\frac{2^{n+2}-\vert\e'\vert}{2^{n+1}}.$$

$$\begin{array}{ccl}
        \Vert f(\e)-f(\e')\Vert & = & \Vert \lambda
        f_{n}(\e)+(1-\lambda)f_{n+1}(\e)-\lambda'f_{n+1}(\e')-(1-\lambda')f_{n+2}(\e')\Vert\\
         & & \\
         & \le & \lambda(\Vert f_{n}(\e)\Vert+\Vert
         f_{n+1}(\e)\Vert)+(1-\lambda')(\Vert f_{n+1}(\e')\Vert+\Vert f_{n+2}(\e')\Vert)+\Vert f_{n+1}(\e)-f_{n+1}(\e')\Vert \\
          & & \\
          & \le &
          \dis+2\lambda\vert\e\vert+2(1-\lambda')\vert\e'\vert\\
           & & \\
           & \le & 9\dis,
         \end{array}$$
because,
$$\lambda \leq \frac{\dis}{2^n},\ \ {\rm so}\ \ \lambda\vert\e\vert\leq 2\dis.$$
Similarly
$$1-\lambda'=\frac{\vert\e'\vert-2^{n+1}}{2^{n+1}}\leq \frac{\dis}{2^{n+1}}\ \ {\rm
and}\ \ (1-\lambda')\vert\e'\vert\leq 2\dis.$$

Finally, $f$ is $9$-Lipschitz.\\

\end{enumerate}

Now we deal with the minoration.

\noindent In our next discussion, whenever $\vert\e\vert$ (respectively
$\vert\e'\vert$) will belong to $[2^{n},2^{n+1})$, for some integer $n$, we
shall denote
$$\lambda=\frac{2^{n+1}-\vert\e\vert}{2^{n}}\ \ ({\rm respectively}\ \ \lambda'=\frac{2^{n+1}-\vert\e'\vert}{2^{n}}).$$
We can suppose that $\e$ is smaller than $\e'$ in the lexicographic order.
Denote $\delta$ the greatest common ancestor of $\e$ and $\e'$. And let
$d=\vert\e\vert-\vert\delta\vert$ (respectively
$d'=\vert\e'\vert-\vert\delta\vert$).

\begin{enumerate}[1)]

\item  if $2^{n}\le\vert\e\vert,\vert\e'\vert\le 2^{n+1}$.\\

We have,

$$x_{n,\Psi_{n}(\delta)}^{*}\Pi_{n}(f(\e)-f(\e'))=\theta(\lambda
d-\lambda' d')$$

$$x_{n+1,\Psi_{n+1}(\delta)}^{*}\Pi_{n+1}(f(\e)-f(\e'))=\theta((1-\lambda)
d-(1-\lambda') d').$$

Hence, $$\Vert f(\e)-f(\e')\Vert\ge\frac{\theta(d-d')}{8}.$$

And,

$$-x_{n,\Psi_{n}(\e)}^{*}\Pi_{n}(f(\e)-f(\e'))=\theta\lambda'
d'$$

$$-x_{n+1,\Psi_{n+1}(\e)}^{*}\Pi_{n+1}(f(\e)-f(\e'))=\theta(1-\lambda')
d'.$$

So, $$\Vert f(\e)-f(\e')\Vert\ge\frac{\theta d'}{8}.$$

Finally if we distinguish the cases $\frac{d}{2}\le d'$, and $d'<\frac{d}{2}$
we obtain :\\

$$\Vert f(\e)-f(\e')\Vert\ge\frac{\theta(d+d')}{24}=\frac{\theta}{24}\dis.$$

\item if $2^{n}\le\vert\e\vert\le 2^{n+1}\le 2^{q+1}\le\vert\e'\vert\le
2^{q+2}$,\\
or $2^{n}\le\vert\e'\vert\le 2^{n+1}\le 2^{q+1}\le\vert\e\vert\le
2^{q+2}$.\\

If $n<q$,

$$\vert x_{q+1,\Psi_{q+1}(\delta)}^{*}\Pi_{q+1}(f(\e)-f(\e'))+x_{q+2,\Psi_{q+2}(\delta)}^{*}\Pi_{q+2}(f(\e)-f(\e'))\vert=\theta Max(d,d')$$

Hence,

$$\Vert f(\e)-f(\e')\Vert\ge\frac{\theta}{16}\dis.$$

If $n=q$ and $\vert\e\vert\le\vert\e'\vert$,

$$\vert
x_{n+1,\Psi_{n+1}(\e)}^{*}\Pi_{n+1}(f(\e)-f(\e'))+x_{n+2,\Psi_{n+2}(\delta)}^{*}\Pi_{n+2}(f(\e)-f(\e'))\vert\ge\theta
d'.$$

So, $$\Vert f(\e)-f(\e')\Vert\ge\frac{\theta}{16}\dis.$$

If $n=q$ and $\vert\e'\vert<\vert\e\vert$,

$$x_{n+1,\Psi_{n+1}(\delta)}^{*}\Pi_{n+1}(f(\e)-f(\e'))-x_{n+1,\Psi_{n+1}(\e)}^{*}\Pi_{n+1}(f(\e)-f(\e'))+x_{n+2,\Psi_{n+2}(\delta)}^{*}\Pi_{n+2}(f(\e)-f(\e'))=\theta
d.$$

Hence,

$$\Vert f(\e)-f(\e')\Vert\ge\frac{\theta}{24}\dis.$$

\end{enumerate}

Finally $T\buildrel {\frac{216}{\theta}}\over {\hookrightarrow} X$.

\end{pf}

\begin{cor}
$X$ is non super-reflexive if and only if $(T,\rho)$ embeds into $X$.
\end{cor}

\begin{pf}
It follows clearly from Bourgain's result \cite{B} and Theorem \ref{theo1}.
\end{pf}
\\

\section{Metric characterization of the linear type}

First we identify canonicaly $\{-1,1\}^{n}$ with
$K_n=\{-1,1\}^{n}\times\prod_{k>n}\{0\}$.\\

Let $p\in[1,\infty)$.\\
Then we define an other metric on $T=\bigcup K_n$ as follows :\\
$\forall\ \e,\e'\in T$,
$$d_{p}(\e,\e')=(\sum_{i=0}^{\infty}\vert\e_{i}-\e'_{i}\vert^p)^{\frac{1}{p}}.$$

The length of $\e\in T$ can be viewed as $\vert\e\vert=(d_{p}(\e,0))^p$.\\
The norm $\Vert . \Vert_{p}$ on $\ell_{p}$ coincides with $d_{p}$ for the
elements of $T$.\\

\noindent We recall now two classical definitions :\\

Let $X$ and $Y$ be two Banach spaces. If $X$ and $Y$ are linearly isomorphic,
the {\it Banach-Mazur distance} between $X$ and $Y$, denoted by $d_{BM}(X,Y)$,
is the infimum of $\|T\|\,\|T^{-1}\|$, over all linear isomorphisms $T$ from
$X$ onto $Y$.\\

For $p\in [1,\infty]$, we say that a Banach space $X$ uniformly contains the
$\ell_p^n$'s if there is a constant $C \geq 1$ such that for every integer $n$,
$X$ admits an $n$-dimensional subspace $Y$ so that $d_{BM}(\ell_p^n,Y)\leq
C$.\\

We state  and prove now the following result.

\begin{thm}\label{lp}
Let $p\in [ 1,\infty ) $.\\
If $\ X$ uniformly contains the $\ell_{p}^{n}$'s then $(T,d_{p})$ embeds into
$X$.\\
\end{thm}

\begin{pf}
We first recall a fundamental result due to Krivine (for $1<p<\infty$ in
\cite{Kr}) and James (for $p=1$ and $\infty$ in \cite{J}).

\begin{thm}[James-Krivine]\label{fcod} Let $p\in[1,\infty]$ and $X$ be a Banach space uniformly containing the $\ell_p^n$'s.
Then, for any finite codimensional subspace $Y$ of $X$, any $\epsilon>0$ and
any $n\in \N$, there exists a subspace $F$ of $Y$ such that
$d_{BM}(\ell_p^n,F)<1+\epsilon$.
\end{thm}

Using Theorem \ref{fcod} together with the fact that each $\ell_{p}^{n}$ is
finite dimensional, we can build inductively finite dimensional subspaces
$(F_{n})_{n=0}^{\infty}$ of $X$ and $(R_{n})_{n=0}^{\infty}$ so that for every
$n\geq 0$, $R_{n}$ is a linear isomorphism from $\ell_{p}^{n}$ onto $F_{n}$
satisfying
$$\forall u\in \ell_{p}^{n}\ \ \ \ \frac{1}{2}\|u\|\leq \|R_{n}u\|\leq \|u\|$$
and also such that $(F_{n})_{n=0}^{\infty}$ is a Schauder finite dimensional
decomposition of its closed linear span $Z$. More precisely, if $P_{n}$ is the
projection from $Z$ onto $F_0\oplus...\oplus F_{n}$ with kernel $\overline{\rm
Span}\,(\bigcup_{i=n+1}^{\infty} F_{i})$, we will assume as we may, that
$\|P_{n}\|\leq 2$. We denote now $\Pi_{0}=P_{0}$ and $\Pi_{n}=P_{n}-P_{n-1}$
for $n\geq 1$. We have that $\|\Pi_{n}\|\leq 4$.

We now consider $\varphi_{n}:T_{n} \to \ell_{p}^{n}$ defined by
$$\forall \e\in T_{n},\ \ \varphi_{n}(\e)=
\sum_{i=1}^{\vert\e\vert}\e_{i}e_{i},$$ where $(e_{i})$ is the canonical basis
of $\ell_{p}^{n}$. The map $\varphi_{n}$ is clearly an isometric embedding of
$T_{n}$ into $\ell_{p}^{n}$.

\noindent Then we set :
$$\forall \e\in T_{n},\ \ f_{n}(\e)=R_{n}(\varphi_{n}(\e)) \in F_{n}.$$

\noindent Finally we construct a map $f:T\to X$ as follows :

$$\begin{array}{cccl}

    f : & T & \rightarrow &  X \\
        &   &   &  \\
        & \e & \mapsto & \lambda f_{m}(\e)+(1-\lambda)f_{m+1}(\e)\ ,\ {\rm if}\ 2^{m}\le\vert\e\vert< 2^{m+1},\\

    \end{array}$$

where, $$\lambda=\frac{2^{m+1}-\vert\e\vert}{2^{m}}.$$\\

\begin{rk}
We have $\frac{1}{16}\vert\e\vert^{\frac{1}{p}}\le\Vert
f(\e)\Vert\le\vert\e\vert^{\frac{1}{p}}$.
\end{rk}
\medskip Like in the proof of Theorem \ref{theo1} ,we prove that $f$ is $9$-Lipschitz using
exactly
the same computations.\\

\medskip We shall now prove that $f^{-1}$ is Lipschitz. We consider
$\e,\e'\in T$ and assume again that $0<\vert\e\vert\leq \vert\e'\vert$. We need
to study two different cases. Again, whenever $\vert\e\vert$ (respectively
$\vert\e'\vert$) will belong to $[2^{m},2^{m+1})$, for some integer $m$, we
shall denote
$$\lambda=\frac{2^{m+1}-\vert\e\vert}{2^{m}}\ \ ({\rm respectively}\ \
\lambda'=\frac{2^{m+1}-\vert\e'\vert}{2^{m}}).$$\\

\begin{enumerate}[1)]

\item if $2^{m}\le\vert\e\vert,\vert\e'\vert<
2^{m+1}$.\\

$$\begin{array}{ccl}
d_{p}(\e,\e') & \le &
\Vert\lambda\sum_{i=1}^{\vert\e\vert}\e_{i}e_{i}-\lambda'\sum_{i=1}^{\vert\e'\vert}\e'_{i}e_{i}\Vert_{p}+\Vert(1-\lambda)\sum_{i=1}^{\vert\e\vert}\e_{i}e_{i}-(1-\lambda')\sum_{i=1}^{\vert\e'\vert}\e'_{i}e_{i}\Vert_{p}\\
 & & \\
 & \le & 2\Vert\Pi_{m}(f(\e)-f(\e'))\Vert+2\Vert\Pi_{m+1}(f(\e)-f(\e'))\Vert\\
 & & \\
  & \le & 16\Vert f(\e)-f(\e')\Vert.\\
\end{array}$$

\item if $2^{m}\le\vert\e\vert\le 2^{m+1}\le 2^{q+1}\le\vert\e'\vert<
2^{q+2}$.\\

\indent if $m<q$,

$$\begin{array}{ccl}
d_{p}(\e,\e') & \le & 2 d_p(\e',0)\\
 & & \\
 & \le & 2((1-\lambda')d_p(\e',0)+\lambda'd_p(\e',0))\\
 & & \\
 & \le & 2(2\Vert\Pi_{q+2}(f(\e)-f(\e'))\Vert+2\Vert\Pi_{m+1}(f(\e)-f(\e'))\Vert)\\
 & & \\
  & \le & 32\Vert f(\e)-f(\e')\Vert.\\
\end{array}$$

\newpage

\indent if $m=q$,

$$\begin{array}{ccl}
d_{p}(\e,\e') & \le & \lambda d_p(\e,0)+
\Vert(1-\lambda)\sum_{i=1}^{\vert\e\vert}\e_{i}e_{i}-\lambda'\sum_{i=1}^{\vert\e'\vert}\e'_{i}e_{i}\Vert_{p}+(1-\lambda')d_p(\e',0)\\
 & & \\
 & \le & 2\Vert\Pi_{m}(f(\e)-f(\e'))\Vert+2\Vert\Pi_{m+1}(f(\e)-f(\e'))\Vert+2\Vert\Pi_{m+2}(f(\e)-f(\e'))\Vert\\
 & & \\
  & \le & 24\Vert f(\e)-f(\e')\Vert.\\
\end{array}$$

\end{enumerate}
Finally we obtain that $f^{-1}$ is $32$-Lipschitz, and $T\buildrel {288}\over
{\hookrightarrow} X$.

\end{pf}

In the sequel a Banach space $X$ is said to have {\it type $p>0$} if there
exists a constant $T<\infty$ such that for every $n$ and every
$x_1,\dots,x_n\in X$,
$$\mathbb{E}_{\e}\Vert\sum_{j=1}^{n}\e_j x_j\Vert^p_X\le T^p\sum_{j=1}^{n}\Vert
x_j\Vert_X^p,$$ where the expectation $\mathbb{E}_\e$ is with respect to a
uniform choice of signs $\e_1,\dots,\e_n\in\{-1,1\}^n$.

The set of $p$'s for which $X$ contains $\ell_{p}^{n}$'s uniformly is closely
related to the type of $X$ according to the following result due to Maurey,
Pisier \cite{MP} and Krivine \cite{Kr}, which clarifies the meaning of these
notions.

\begin{thm}[Maurey-Pisier-Krivine]
Let $X$ be  an infinite-dimensional Banach space. Let $$p_X={\rm sup}\{p\ ;{\rm
X\ is\ of\ type\ p\}},$$\\
Then $X$ contains $\ell_{p}^{n}$'s uniformly for $p=p_X$.\\
Equivalently, we have $$p_X={\rm inf}\{p\ ;{\rm X\ contains\
\ell_{p}^{n}\mbox{'s}\ uniformly\}}.$$

\end{thm}

We deduce from Theorem \ref{lp} two corollaries.

\begin{cor}
Let $X$ a Banach space and $1\le p<2$.\\
The following assertions are equivalent :

\begin{enumerate}[i)]
\item $p_X\le p$.
\item $X$ uniformly contains the $\ell_{p}^{n}$'s.
\item the $(T_{n},d_{p})$'s uniformly embed into $X$ .
\item $(T,d_{p})$ embeds into $X$.
\end{enumerate}

\end{cor}

\begin{pf}
$ii)$ implies $i)$ is obvious.\\

$i)$ implies $ii)$ is due to Theorem \ref{fcod} and the work of Bretagnolle, Dacunha-Castelle and Krivine \cite{BDCK}.\\

For the equivalence between $ii)$ and $iii)$ see the work of Bourgain, Milman and Wolfson \cite{BMW} and Krivine \cite{Kr}.\\

$iv)$ implies $iii)$ is obvious.\\

And $ii)$ implies $iv)$ is Theorem \ref{lp}.

\end{pf}

\begin{cor}
Let $X$ be an infinite dimensional Banach space, then $(T,d_2)$ embeds into
$X$.
\end{cor}

\begin{pf}
This corollary is a consequence of the Dvoretsky's Theorem \cite{Dv} and
Theorem \ref{lp}.

\end{pf}


\begin{thebibliography}{99}
\bibitem{A} I. Aharoni, {\it Every separable metric space is Lipschitz equivalent to a subset of $c_{0}^+$}. Israel J. Math. 19 (1974),
284--291.\\

\bibitem{B} J. Bourgain, {\it The metrical interpretation of super-reflexivity in Banach spaces}. Israel J. Math.  56 (1986),
221-230.\\

\bibitem{BMW} J. Bourgain, V. Milman, H. Wolfson, {\it On type of metric
spaces}. Trans. Amer. Math. Soc. volume 294, number
1, march 1986, 295-317.\\

\bibitem{BDCK} J. Bretagnolle, D. Dacunha-Castelle, J.L. Krivine,
{\it Lois stables et espaces $L^p$}. Ann. Instit. H. Poincar\'{e}, 2
(1966), 231-259.\\

\bibitem{D} J. Diestel, Sequences and Series in Banach Spaces. Springer-Verlag
(1984).\\

\bibitem{Dv} A. Dvoretzky, {\it Some results on convex bodies and Banach spaces}. Proc. Internat. Sympos. Linear Spaces (Jerusalem, 1960)
123--160.\\

\bibitem{J} R. C. James, {\it Super-reflexive spaces with bases}. Pacific J. Math. 41
(1972), 409-419.\\

\bibitem{Kr} J. L. Krivine, {\it Sous-espaces de dimension finie des
espaces de Banach r\'{e}ticul\'{e}s}, Ann. of Math. (2) 104 (1976), 1-29.\\

\bibitem{LT} J. Lindenstrauss and L. Tzafriri, Classical Banach Spaces I,
Springer Berlin 1977.\\

\bibitem{MP} B. Maurey, G. Pisier, {\it S\'{e}ries de variables al\'{e}atoires vectorielles
ind\'{e}pendantes et propri\'{e}t\'{e}s g\'{e}om\'{e}triques des espaces de Banach},
Studia Math. {\bf 58(1)} 1976, 45-90.\\

\bibitem{MN} M. Mendel, A. Naor, {\it Metric cotype}, arXiv:math.FA/0506201 v3 29 Apr 2006\\

\bibitem{Pi} G. Pisier, Factorization of Linear Operators and
Geometry of Banach Spaces. CBMS Regional Conference Series in Mathematics,
60.\\

\bibitem{R} M. Ribe, {\it Existence of separable uniformly homeomorphic non
isomorphic Banach spaces}.  Israel J. Math. 48 (1984), no. 2-3, 139-147.\\

\end{thebibliography}
\end{document}